\documentclass[twoside,12pt]{article}
\usepackage{amsmath,amssymb,epsfig}
\date{}
\newtheorem{Lemma}{LEMMA}[section]
\newtheorem{Corollary}[Lemma]{COROLLARY}
\newtheorem{Theorem}[Lemma]{THEOREM}
\newtheorem{Proposition}[Lemma]{PROPOSITION}

\newcommand{\bnum}{\begin{enumerate}}
\newcommand{\enum}{\end{enumerate}}
\newcommand{\bi}{\begin{itemize}}
\newcommand{\ei}{\end{itemize}}
\newcommand{\btab}{\begin{tabular}}
\newcommand{\etab}{\end{tabular}}
\newcommand{\beq}{\begin{eqnarray*}}
\newcommand{\eeq}{\end{eqnarray*}}
\newcommand{\beqn}{\begin{eqnarray}}
\newcommand{\eeqn}{\end{eqnarray}}

\newcommand{\bq}{\begin{equation}}
\newcommand{\eq}{\end{equation}}

\newcommand{\CF}{{\cal F}}
\newcommand{\CH}{{\cal H}}

\newcommand{\CL}{{\cal L}}

\newcommand{\CP}{{\cal P}}

\newcommand{\CS}{{\cal S}}

\def\phi{\varphi}
\def\epsilon{\varepsilon}
\newcommand{\BS}{\Bbb S}
\newcommand{\BR}{\Bbb R}
\newcommand{\BC}{\Bbb C}

\newcommand{\BT}{\Bbb T}
\newcommand{\Pd}{\mathrm{PG}(3,\BR)}
\newcommand{\Pf}{\mathrm{PG}(5,\BR)}
\newcommand{\mo}{\mathop\mathrm}
\newcommand{\kasten}{\vbox{\hrule height 8pt width 8.6pt depth -7.4pt
    \hbox{\vrule width 0.6pt height 7.4pt
    \kern 7.4pt \vrule width 0.6pt height 7.4pt}
    \hrule height 0.6pt width 8.6pt}}
\newcommand{\ok}{\hfill\kasten}
\newcommand{\bpf}{\begin{Proof}}
\newcommand{\epf}{\ok\end{Proof}\bigskip\noindent}
\newcommand{\bthm}{\begin{Theorem}}
\newcommand{\ethm}{\end{Theorem}}
\newcommand{\ble}{\begin{Lemma}}
\newcommand{\ele}{\end{Lemma}}
\newcommand{\bprop}{\begin{Proposition}}
\newcommand{\eprop}{\end{Proposition}}
\newcommand{\bcor}{\begin{Corollary}}
\newcommand{\ecor}{\end{Corollary}}
\begin{document}
\title{Regular parallelisms on $\Pd$ admitting a 2-torus action}
\author{Rainer L\"owen and G\"unter F. Steinke}

\maketitle
\thispagestyle{empty}

\begin{abstract}
A regular parallelism of real projective 3-space $\Pd$ is an equivalence relation on the line space such that every class 
is equivalent to the set of 1-dimensional complex subspaces of $\BC^2 = \BR^4$. We shall assume that the set of classes is 
compact, and characterize those regular parallelisms that admit an action of a 2-dimensional torus group. 
We prove that there is a one-dimensional subtorus fixing every parallel class. From this property alone 
we deduce that the parallelism is
a 2- or 3-dimensional regular parallelism in the sense of Betten and Riesinger \cite{Stevin}. If a 2-torus acts, then the parallelism
can be described using a so-called generalized line star or \it gl star \rm which admits a 1-torus action. 
We also study examples of such parallelisms
by constructing gl stars. In particular, we prove a claim which was presented  in \cite{Stevin} with an incorrect proof.
The present article continues a series of papers by the first author on parallelisms with large groups.
\\
 
\end{abstract}

\section{Introduction}\label{intro}

A \it spread \rm in real projective 3-space $\Pd$ is a set of mutually disjoint lines covering the point set.
A \it parallelism \rm on  $\Pd$ is a set $\Pi$ of mutually disjoint spreads covering the set of lines. A parallelism
is usually considered as an equivalence relation on the line set. We shall always assume
that spreads and parallelisms are topological in the sense explained in Section \ref{basic}. 
That property is equivalent to compactness 
in some natural topology. The classical Clifford parallelism consists of the orbits of one of the two factors of 
$\mo{SO}(4,\BR) = \mo{SO}(3,\BR)\cdot \mo{SO}(3,\BR)$.
The pioneering work of Betten and Riesinger (e.g., \cite{Thaswalker}) 
has produced numerous ways of constructing nonclassical examples.

The automorphism group $\mathop \mathrm{Aut}\Pi \le \mathop \mathrm{PGL}(4,\BR)$
of a topological parallelism is a compact Lie group \cite{beloe}, \cite {comp full gp} of dimension at most 6. There are plenty 
examples with a small group (of dimension $\le 1$), see \cite{bsp aut pi}, 
and there is no hope of knowing them all. At the other end, if the group
is of dimension $> 3$, then the parallelism is Clifford, see \cite{charcliff}, and its group is 6-dimensional. 
Parallelisms with a 3-dimensional group are completely 
known \cite{so3}, \cite{su2}.  Passing from dimension 3 to 2, matters suddenly become extremely difficult. 
So it seems reasonable at first to restrict 
attention to a particularly nice type of parallelisms, the regular ones, where every parallel class is a regular spread, i.e. 
isomorphic to the 
set of complex one-dimensional subspaces of $\BC^2$. This has the advantage that these parallelisms 
can be handeled conveniently using the 
Klein quadric $K$ in $\Pf$, where regular spreads appear as quadrics obtained by intersecting  $K$ with 
3-dimensional subspaces of $\Pf$. Moreover, the  group of automorphisms and dualities of $\Pd$, 
which contains the group of the parallelism,
corresponds to the group $\mathop \mathrm{PO} (6,3)$, compare \cite{Stevin}, and the 2-tori inside this group are easy to see. 

In order to specify  a regular parallelism  $\Pi$, one has to specify the set of 3-spaces used to define 
quadrics in $K$. Equivalently,
one can specify the polar lines of those 3-spaces with respect to the polarity $\pi_5$ which defines the Klein quadric. 
These lines are 0-secants of $K$, i.e., disjoint from $K$. A set $\cal H$ of 
0-secants of $K$ defines 
a parallelism in this way if and only if it is a so called hfd line set, i.e., if every tangent hyperplane $\pi_5(p)$, $p\in K$, 
contains exactly one line
in $\cal H$, see \cite{hyper} or \cite{gldirect}. 
The parallelism is topological if and only if the set $\cal H$ is compact. Betten and Riesinger 
introduced the term \it dimension \rm  of $\Pi$ to mean the dimension of the subspace of $\Pf$ generated by the hfd line 
set $\cal H$ corresponding to $\Pi$,
$$\dim \Pi = \dim\mathop\mathrm{span} \cal H.$$

Examples exist for 
parallelisms of dimensions 2, 3, 4 and 5, see \cite{latitud}, \cite{hyper}. Dimension 2 charcterizes Clifford 
parallelism; in this case, $\cal H$ consists of 
\it all \rm lines of the plane generated by $\cal H$. Our main result is

\bthm\label{main} \rm
A regular topological parallelism with a 2-dimensional automorphism group is at most 3-dimensional.
\ethm

The only 2-dimensional compact connected Lie group is the 2-torus. Examples of regular parallelisms 
with 2-torus action are known only in very special cases. 
We shall give simplified proofs of the existing results in Sections \ref{axi} and \ref {sym case}, and we shall prove strong 
general existence results in Section \ref {general}.

In the 3-dimensi\-onal case, the hfd line set $\cal H$ 
defining a regular parallelism can be replaced by a simpler object, namely, 
by a compact so-called \it generalized line star \rm or \it gl star $\cal S$. \rm 
The gl star is obtained by 
applying to the set $\cal H$ of lines  the polarity $\pi_3$ induced by $\pi_5$ on 
the 3-space $P_3$ generated by $\cal H$,
$${\cal S} = \pi_3 \pi_5(\Pi).$$
This is due to Betten and Riesinger \cite{gl}; see \cite{gldirect} for a simpler approach and proof. The defining property 
of a gl star, which is equivalent to the fact that the gl star corresponds to a 3-dimensional parallelism, is this: $\cal S$ is a set 
of lines in a 3-dimensional subspace $P_3$ of $\Pf$, and $Q = P_3 \cap K$ is an elliptic quadric such that every line in $\cal S$
meets $Q$ in two points, and every point of $P_3$ not in the interior of $Q$ is incident with precisely one line from $\cal S$.
Here, a point is called exterior if it is incident with an exterior line, that is, a 0-secant of $Q$. An interior point is a point 
$p$ such that every line containing $p$ is a 2-secant (meeting $Q$ in two points). The remaining points are precisely those of $Q$.

We shall show (Theorem \ref{crit}) that for compact $\cal S$ the defining property of gl-stars can be simplified: 
it suffices that every point on $Q$ is on some line of $\cal S$ and two lines from $\cal S$ 
never meet in a non-interior point with respect to $Q$. 

A compact gl star can be described by a  fixed point free 
involutory homeomorphism $\sigma$ of $Q$; the star consists of the 
lines $q \vee \sigma(q)$ for all $q\in Q$, and conversely, $\sigma$ sends
a point $q \in Q$ to the second point of intersection of $Q$ with the line $L \in \cal S$ containing $q$. The gl star 
corresponding to Clifford parallelism consists of all lines passing to some fixed point $p$ in the interior of $Q$, 
it is an \it ordinary line star. \rm The corresponding involution is the antipodal map with respect to $p$.


\section{Basic facts and known results}\label{basic}

We shall freely use the continuity properties of the topological projective space $\Pd$, see, e.g., \cite{Kuehne}.

We recall the definition and the properties of the \it  Klein correspondence \rm in the case of real projective 3-space $\Pd$,
following Pickert \cite{Pickert}; compare also \cite{Knarr}, Section 2.1. 
By definition, $\Pd$ is the lattice of subspaces of $\BR^4$. In the 6-dimensional vector space of non-degenerate 
alternating forms on $\BR^4$, we can represent a line $L$ of $\Pd$ by the one-dimensional subspace of forms containing $L$ in their radical.
This maps the line set $\CL$ of $\Pd$ bijectively onto a quadric $K$ in $\Pf$, called the \it Klein quadric, 
\rm defined by a nondegenerate quadratic form $g$ of index 3.
The space $\BR^6$ contains two families of maximal isotropic subspaces (of projective dimension 2), 
representing the points and the hyperplanes of $\Pd$, respectively.
Incidence between lines and points or hyperplanes in $\Pd$ corresponds to inclusion in the Klein model of $\Pd$.

\it Automorphisms \rm of $\Pd$ correspond to orthogonal maps with respect to $g$. For us, it suffices to know that the connected
component $\mathop \mathrm{PSL} (\BR,4)$ of $\mathop \mathrm {Aut}\Pd$ corresponds to the projective orthogonal group 
$\mo {PSO}(6,3)$ with respect to $g$. The automorphism group $\mathop \mathrm {Aut}\Pi$ of a regular parallelism appears 
in the Klein model as a subgroup of this group. According to \cite{beloe}, \cite{comp full gp}, this group is compact 
and hence is contained in the maximal compact subgroup $\mathrm {O}(3,\BR) \cdot \mathrm {O}(3,\BR)$ (product with 
amalgamated central involutions).
The factors of this product are induced by the orthogonal groups on the factors of $\BR^6 = \BR^3 \times \BR^3$, 
where the quadratic form defining the Klein quadric is given by $g(u,v) = \Vert u \Vert^2 - \Vert v \Vert^2$, with 
$\Vert u \Vert$ denoting the Euclidean norm. A 2-torus in $\mo{Aut} \Pi$ is therefore a product $T_1 \cdot T_2$, 
where $T_i\cong \mo{SO}(2,\BR)$ consists of the rotations of the $i$th factor $\BR^3$ about some fixed axis.
Clearly, the $g$-orthogonal group $\mo{Aut}\Pi$ leaves the 3-space $P_3 = \mo{span} \cal H$ and the gl star $\cal S$ 
associated with $\Pi$ invariant.

By definition, a spread of $\Pd$ is \it regular \rm if together with any triple of distinct lines it contains the entire regulus 
determined by this triple.
As pointed out earlier, among spreads in $\Pd$ this property charcterizes the complex spread, 
consisting of the one-dimensional complex subspaces with respect to some complex structure on $\BR^4$, see \cite{Knarr}, 4.15.
A spread of $\Pd$ is regular if and only if its
image under the Klein correspondence is an elliptic subquadric of $K$ spanning a 3-space, see, e.g., \cite{Stevin}, Proposition 13. 

\it Topological \rm spreads, topological regular parallelisms, topological gl stars and topological hfd line sets are discussed 
thoroughly in \cite {gldirect}, Section 3; here
we give a brief summary. In the case of a parallelism, being topological means that the unique line parallel (i.e., equivalent) 
to a given line $L$ and passing through a given point $p$ depends continuously on the pair $(p,L)$. For a regular parallelism, 
this is equivalent to compactness of the set of 3-spaces in $\Pf$ defining the regular spreads belonging to the parallelism. 
The topology to be considered on this set is the topology of the Grassmann manifolds of all 3-spaces in $\Pf$.
For the other objects enumerated above, the definition of a topological object is similar; e.g., for a spread it is 
required that the line containing a given point depends continuously on that point. Again, compactness in the topology 
induced by the relevant Grassmann manifold is equivalent to this property.\\

About the possibilities for automophism groups of regular parallelisms, the following results are known. We begin with large groups.

\bthm\label{dim ge 3}
A regular parallelism with an automorphism group of dimension at least 3 is Clifford.
\ethm

\bpf
This is a corollary of the first author's results on general parallelisms with a 3-dimensional group \cite{su2}, Corollary 3.2. 
For 3-dimensional regular parallelisms, a simple proof is available, compare \cite{Stevin}, Theorem 33. At the time, 
the compactness of automophism groups was not known, therefore we update the argument here: A compact subgroup of $\mo{PO}
(6,3)$ of dimension at least 3 contains 
a copy of $\mo{SO}(3,\BR)$, which acts on the quadric $Q \cong \BS_2$ associated with the gl star $\cal S$ in the ordinary way. 
The rotations fixing a point $q \in Q$ fix the unique line $L = q \vee \sigma (q)\in \CS$ that contains $q$, hence $L$ coincides 
with the rotation axis 
and passes through the universal fixed point of the rotation group. Thus we have an ordinary gl star and a Clifford parallelism. 
\epf

On the other hand, a 3-dimensional regular parallelism always has an automorphism group of positive dimension:

\bthm\label{ex 1-torus}
\cite{Stevin} Every at most 3-dimensional regular parallelism admits an action of a 1-torus $\BT$ that fixes every parallel class. 
\ethm

Actually, the same proof yields an action of the orthogonal group $\mathrm{O}(2,\BR)$ with the same property.

\bpf
We may assume that the given regular parallelism $\Pi$ is not Clifford. Hence $P_3 = \mo{span} \CH$ is a 3-space, and on the 
underlying vector space $U = \BR^4$ the form $g$ induces a form of signature $(3,1)$. Thus the $g$-orthogonal complement $C$ 
of $U$
in $\BR^6$ is of signature $(0,2)$. There is a $g$-orthogonal effective torus action that 
induces $\mo{SO}(2,\BR)$ on $C$   and is trivial
on $U$. Therefore it fixes all elements of $\CS$ and, hence,  all parallel classes.
\epf

In Section \ref{2-torus}, 
we shall prove the converse of this theorem, and this will be the key to the proof of our main result \ref{main}.\\

Betten and Riesinger have given numerous \it examples \rm of non-Clifford regular parallelisms and their automorphism groups. 
In \cite{latitud} they construct 3-dimensional regular parallelisms admitting a 2-torus. 
These examples are of a very special kind called axial.
A much simpler treatment of these examples and their properties will be given in Section \ref{axi}. 
In \cite {Stevin}, Theorem 42, they give a 
construction which they claim yields more general examples with 2-torus action. Their proof is invalid, 
but we shall prove an improved version of their claim (Theorem \ref{42}).
In \cite{bsp aut pi}, Betten and Riesinger present regular parallelisms of 
dimension 3 and 4 whose full automorphism group is 1-dimensional, 
as well as regular parallelisms of dimension 4 and 5 whose  group is 0-dimensional (hence finite). 
In \cite{bsp aut pi} they ask whether groups of dimension 2 can act on regular parallelisms of dimension 4 or 5. Our main 
result Theorem \ref {main} answers this question in the negative. They also ask the following\\

\bf Question: \rm Can a 1-torus act on some 5-dimensional regular parallelism? \\

This question remains open.


\section{2-Torus action: The proof of Theorem \ref{main}}\label{2-torus}

We need the following well-known fact.

\bprop \label{ineff}
A 2-torus $\BT^2$ cannot act effectively on a connected surface other than a torus.
\eprop

\bpf
Let $S$ be the surface and assume first that $\BT^2$ has a 2-dimensional orbit $B$. Then $B$ is homeomorphic to a torus and is 
open and closed in the connected surface $S$, hence $S = B$ is a torus. If all orbits of $\BT^2$ are of dimension at most 1, we 
can argue in two independent ways. First, all stabilizers are of dimension at least one, and there are only countably many such 
subgroups in the torus. By Baire's theorem, one of them has a fixed point set with nonempty interior. By Newman's 
theorem \cite {newm}, \cite {dress}, this subgroup acts trivially. 

The other argument uses the results of Mostert \cite {mos} on 
compact group actions with orbits of codimension one. The orbits with minimal stabilizers (called principal orbits) cover a dense open set,
and all stabilizers of points on principal orbits are conjugate, hence in fact equal. This subgroup then acts trivially on $S$.
\epf

Note that a 2-dimensional factor group of a 2-torus is again a 2-torus, hence the kernel of ineffectivity has positive 
dimension if a 2-torus acts on a connected surface which is not a torus. A parallelism $\Pi$ is homeomorphic to a star of 
lines in $\Pd$ via the map that sends a line to its unique parallel passing through 
some fixed point. Thus $\Pi$ is homeomorphic to the real projective plane and certainly not to a torus. It follows that a 2-torus 
cannot act effectively on $\Pi$, and our Theorem \ref{main} becomes a corollary to the  converse of Theorem \ref{ex 1-torus}, 
which we prove next.

\bthm \label {trivial torus}
Let $\Pi$ be a topological regular parallelism. Then $\mo{Aut}\Pi$ contains a 1-torus $\BT$ that fixes all parallel classes if and only if 
$\Pi$ is at most 3-dimensional.
\ethm

\bpf 
From the description of tori given in Section \ref{basic} we see that the action of $\BT$ leaves a decomposition 
$\BR^6 = \BR^3 \times \BR^3$ invariant. On both factors, $\BT$ induces (possibly trivial) rotations about some axis. So we may 
regroup factors and obtain the action of $\BT$ on $\BR^6 = \BC^3 = \BC \times \BC \times \BC$ consisting of the maps
     $$\tau_t: (u,v,w) \to (u, e^{kt}v,e^{lt}w),$$
where $t \in [0,2\pi]$ and $k,l$ are relatively prime integers (possibly zero). We claim that $0 \in \{k,l\}$. 
We assume that this is not the case and aim for a contradiction.

The fact that $\BT$ fixes all parallel classes means that $\BT$ acts trivially on the hfd line set $\CH$ associated to $\Pi$. 
In other words, $\CH$  is contained in the set $\CF$ of all fixed lines of $\BT$ in $\Pf$. In $\BC^3$, these lines appear 
as $\BT$-invariant 2-dimensional real vector subspaces.  So we wish to determine all such subspaces $L$.

If the action of $\BT$ on $L$ is trivial, then $L$ is the first factor of $\BC^3$. So assume that this action is nontrivial. 
Then the $\BT$-equivariant projection of $L$ to the first factor $\BC$ is not bijective, because the actions are different. 
On the other hand, the kernel of this projection is $\BT$-invariant, so the projection must be zero, and $L$ is contained in 
the product $\BC \times \BC$ of the second and third factor. These factors are invariant. If there is another invariant 2-dimensional 
subspace $L$, then by combining the projections of $L$ onto the factors of $\BC \times \BC$ we obtain an equivariant bijection 
between these factors, which means that $k = \pm l$. Since $\BT$ is effective on $\BR^6$, we have $\{k,l\} \subseteq \{1,-1\}$. 
Applying complex conjugation to the factors as necessary, we may obtain $k = l = 1$. Then $\BT$ acts on the complex vector space 
$\BC^2$ via multiplication by complex scalars, and $\CF$ consists of the first factor of $\BC^3$ plus the elements of the complex
(regular) spread of $\BC^2$. Now this spread is homeomorphic to the 2-sphere (the Riemann number sphere). On the other hand, 
$\CH$ is homeomorphic to $\Pi$, i.e., to the real projective plane. Hence the inclusion $\CH \subseteq \CF$ is impossible.

The only remaining possibility is that $k = 0$ and $l = \pm 1$ (up to exchange of factors). 
By considering projections as before we see that then $\CF$ consists of the third factor together with 
all 2-dimensional real subspaces of $\BC \times \BC \times 0$. Since $\CH$ is connected, the third factor does not belong to $\CH$,
and $\mo{span} \CH$ is contained in $\BC \times \BC \times 0$ and has projective dimension at most 3.
\epf

The proof of Theorem \ref{main} is now complete.


\section{Rotational gl stars}\label{rot gl}

If $\Pi$ is a non-Clifford 
regular parallelism admitting an action of a 2-torus $\BT^2$, then we know from Proposition \ref{ineff} that 
some 1-dimensional subtorus acts trivially on $\Pi$, and likewise on the associated gl star $\CS$ and on the hfd line set $\CH$. 
From the proof of Theorem \ref{trivial torus} we infer that it even acts trivially on the 4-dimensional vector space $U$ corresponding 
to $\mo{span} \CH$, on which $g$ has signature $(3,1)$. Up to $g$-orthogonal maps, all such subspaces of $\BR^6$ are 
equivalent by a theorem of Witt, see, e.g.,  \cite{Huppert}, 6.1. 
The entire group $\BT^2$ cannot act trivially on the space $U$, because it cannot act effectively on its 
2-dimensional $g$-orthogonal space in $\BR^6$, which is in fact a vector space complement. 

Thus there is a 1-torus $\Phi$ contained in $\BT^2$ 
which acts nontrivially on the 3-space 
$\mo{span} \CH$ and commutes with the polarity $\pi_3$ induced on it by $\pi_5$. Moreover, $\Phi$ leaves the  quadric $Q$ and the 
gl star $\CS$ invariant. Following Betten and Riesinger \cite{Stevin}, in this situation we call  $\CS$ a \it rotational gl star\rm . 
Up to conjugacy, 
there is only one 1-torus on $\mo{span} \CH$
which leaves $\pi_3$ invariant, so we can assume the following \\

\it Standard position for rotational gl stars: \rm $Q$ is the unit sphere $\BS_2$ in the affine part $\BR^3$ of $\Pd$. 
The group $\Phi$ is the group of 
rotations of $\BR^3$ about the z-axis $Z$, and leaves the gl star $\CS$ invariant. \\

We summarize:

\bprop \label{2-torus-rotational}
\cite{Stevin} Let $\Pi$ be a non-Clifford regular parallelism. Then $\dim \mo{Aut} \Pi = 2$ if 
and only if $\Pi$ is defined by a rotational 
gl star, which we shall always assume in standard position.  \hfill \ok
\eprop

Of course, Clifford parallelism corresponds to a rotational gl star, as well. Note also that $\Phi$ acts on $\CS$ if and 
only if $\Phi$ commutes with the involution $\sigma: \BS_2 \to \BS_2$ defining $\CS$. This will be used in the next proof.

\bprop \label{struct rotational} 
\cite{Stevin} Let $\CS$ be a rotational gl star in standard position.\\
a) The rotation axis $Z$ belongs to $\CS$. \\
b) There is a point $p \in Z$ in the interior of $\BS_2$ such that all lines perpendicular to $Z$ and passing through 
$p$ belong to $\CS$.\\
c) Up to isomorphism, we may assume that the point $p$ is the origin.
\eprop

\bpf
a) The defining involution $\sigma$ commutes with $\Phi$, hence it interchanges the two fixed points $n = (0,0,1)$ and $s = -n$ of 
$\Phi$. Thus $Z = n \vee s \in \CS$.\\
b) For the same reason, $\sigma$ acts on the space $\CS/\Phi \approx [0,1]$ of orbits, interchanging the end points.
Therefore, $\sigma$ sends some $\Phi$-orbit to itself. This implies (b). \\
c) Recall the polarity $\pi_3$ defining the quadric $Q= \BS_2$. The orthogonal group $\Psi = \mo{PSO}(3,1)$ with respect to 
$\pi_3$ is triply transitive on $\BS_2$. The stabilizer $\Psi_{n,s}$  is transitive on the interior open 
segment on $Z$ defined by these points. Moreover, the stabilizer is abelian (isomorphic to $\BC^\times$) and contains $\Phi$. 
Hence we can use it to move the orbit constructed in (b) into the plane $\{(x,y,0)\mid x,y \in \BR\}$ without disturbing 
the action of $\Phi$. 
\epf

In \cite{gl}, a gl star $\CS$ is called \it axial \rm if there is a line $A$, called \it axis, \rm which meets all lines in $\CS$.
We observe:

\bprop\label{axis}
If a non-ordinary (topological) rotational gl star is axial, then the axis as defined above coincides with the rotation axis $Z$.
\eprop

\bpf 
If $A\notin \CS$, then the map $h: \CS \to A$ sending $L\in \CS$ to $L\wedge A$ is continuous, and its image contains the open
subset $B\subseteq A$ of all exterior points on $A$. On $B$, there is a continuous inverse map $j: B \to h^{-1}(B)$, because the 
gl star $\CS$ is topological. Then the surface $\CS$ contains an open subset homeomorphic to the real line, a contradiction. 
Thus $A$ belongs to $\CS$ and meets the rotation axis $Z \in \CS$. Every line $\phi(A)$, $\phi \in \Phi$, has the same properties. 
If $A\ne Z$, then these lines cover either a cone or a plane. Now let $L \in \CS$ not be contained in a plane containing $\Phi A$.
Then $L$ meets all lines $\phi(A)$, and hence must pass through their only common point
$A \wedge Z$, which lies in the interior of $\BS_2$. By density, the remaining lines of $\CS$ also pass through this point, and
$\CS$ is an ordinary line star.
\epf

Rotational gl stars such that every line meets the rotation axis have been called \it latitudinal \rm \cite{latitud}. 
Justified by the above proposition, we shall often use the term \it axial \rm instead. \\

\bthm \label{char ax}
A rotational gl star $\CS$ is axial if and only if the reflection $\zeta_E$ about any plane $E$ containing the axis 
$Z$ acts on $\CS$.  
\ethm

\bpf
The reflection $\zeta_E$ fixes all lines contained in $E$. For $\phi \in \Phi$, the line $\phi(L)$ contained in $\phi(E)$ is 
mapped to $\zeta_E \phi(L) = \zeta_E \phi \zeta_E(L) = \phi^{-1}(L) \in \CS$. Conversely, assume that $\CS = \zeta_E(\CS)$ 
is not axial. Then there is a line 
$L \in \CS$ not meeting $Z$. Choose $\phi \in \Phi$ such that  $\phi(L)$ is parallel to $E$ in the 
Euclidean sense. Then $\phi(L)$ and 
$\zeta_E\phi(L)$ meet at infinity, a contradiction since both lines belong to $\CS$.
\epf


\section{The axial case}\label{axi}

First we simplify the task of recognizing a gl star.

\bprop\label{crit}
Let $\sigma: \BS_2\to \BS_2$ be a continuous fixed point free involution. The set $\CS = \{p \vee \sigma(p) \mid p \in \BS_2\}$ 
of lines in $\Pd$ is a 
topological gl star if  two distinct lines of $\CS$ never meet in a point that is 
exterior with respect to $\BS_2$. In fact, it suffices to assume that every line in $\CS$ 
has a neighborhood in which this condition holds.
\eprop

\bpf
Clearly, $\CS$ is compact. The quotient map $\BS_2 \to \BS_2/\sigma$ is a two sheeted covering map. Hence the quotient space is 
a compact connected surface (in fact, a real projective plane). First suppose that the condition on intersections holds globally.
Then the only property that remains to be verified is that every exterior point $p$ lies on some line of $\CS$.
Let $E$ be a plane containing $p$, disjoint from $\BS_2$. The map 
  $$i_E: \BS_2 \to E: p \to (p\vee \sigma(p))\wedge E$$
induces a continuous injection  $\BS_2/\sigma \to E$ of compact connected surfaces. Such a  map is surjective by domain invariance.
If the condition on intersections holds only locally, then the  map $i_E$ is locally injective, hence a covering map. 
Every covering map between two copies of the real projective plane is bijective, hence the global intersection property follows.
\epf

Now we turn to rotational gl stars $\CS$ that are axial. Let $E$ be a plane containing the rotation axis $Z$. Then  
$\CS$ induces in $E$ a \it gl pencil \rm $\CP$, the 2-dimensional analog of a gl star. That is, $\CP$ covers the plane $E$, and two 
lines of $\CP$ never meet except in the interior of the circle $\BS_2 \cap E$. 
Moreover, $\CP$ is \it symmetric, \rm i.e., invariant under orthogonal reflection $\rho$ in the line $Z$.\\
As for rotational gl stars, we shall mostly assume a \it standard position \rm for symmetric gl pencils: 
the affine part of $E$ is $\BR^2 = \{(x,z)\mid x,z \in \BR\}$, the circle is the unit circle, 
and (as in Proposition \ref{struct rotational}) the coordinate axes $X$ and $Z$ belong to the pencil.\\

In order to give a complete description of all symmetric gl pencils, we introduce the following notation. The intersection of 
the unit circle with 
the first quadrant, i.e., the short segment from $(1,0)$ to $(0,1)$, will be denoted $A$. The opposite segment $-A$ joins the points
$(-1,0)$ and $(0,-1)$. We have the following

\bthm\label{pencil}
a) Let $A \subseteq \BS_1$ be the segment defined above, and let $\mu : A \to -A$ be a continuous bijection which 
sends $(1,0)$ to $(-1,0)$. Then $\mu$ uniquely extends to a fixed point free continuous 
involution $\sigma: \BS_1 \to \BS_1$ that commutes with $\rho$, 
the reflection in the the $z$-axis $Z$. The set
 $$\CP = \{p \vee \sigma(p) \mid p \in \BS_1\}$$
is a topological symmetric gl pencil. \\
b) Every symmetric gl pencil in standard position arises in this way. In particular, every symmetric gl pencil is topological. 
\ethm

\bf Remark: \rm A slightly adapted version of this result holds for non-symmetric gl pencils, with virtually the same proof. We 
have chosen this version because symmetric gl pencils arise from gl stars that are axial and rotational.

\bpf
a) If $p,q$ are two points on the circle, then the pair $(p,\sigma(p))$ separates the pair $(q, \sigma(q))$. 
Indeed, otherwise some segment $B$ of $\BS_1$ bounded by $p$ and $\sigma (p)$ contains $q$ and $\sigma(q)$, hence $B$ is 
mapped to itself and contains a fixed point of $\sigma$, a contradiction.
Therefore, 
any two lines of $\CP$ meet in the interior of the circle and not in the exterior. The proof that every exterior point is covered 
by $\CP$ is virtullay the same as in Proposition \ref{crit}.\\
b) Conversely, if $\CP$ is a symmetric gl pencil (not assumed to be topological) then the separation property 
encountered in (a) holds for $\CP$, because otherwise some pair of $\CP$-lines would meet in the exterior of $\BS_1$. 
This means that the bijection
$\mu : A \to -A$ sending a point $p$ to the other intersection point of $\BS_1$ with the $\CP$-line 
containing $p$ is strictly monotone
with respect to a natural  ordering on $\BS_1 \setminus \{\frac 1 {\sqrt 2}(1,-1)\} \approx \BR$. 
Hence $\mu$ is a homeomorphism and extends 
uniquely to a fixed point free continuous involution $\sigma$ on $\BS_1$ that commutes with $\rho$. 
Then $\CP$ is the gl pencil defined by $\sigma$, and $\CP$ is  compact, hence topological.
\epf

Instead of using the involution $\sigma$, the authors of \cite{latitud} describe the lines of a gl pencil by their intersections 
with the segment $A$ and with the $z$-axis. The pencils are thus given by a function $f$ which expresses the $z$-coordinate of 
the latter point in terms of the $z$-coordinate of the former. This results in rather tricky conditions that $f$ has to 
satisfy and makes it hard to construct examples. In our description, no such conditions appear. \\

Now we can describe all latitudinal (i.e., axial and rotational) gl stars in standard position.

\bthm \label{latitud}
\cite{latitud} a) Let $\CS$ be a latitudinal gl star in standard position and let $E$ be a plane containing the rotation axis $Z$. 
In $E$ then $\CS$ induces  a symmetric gl pencil $\CP$, from which $\CS$ can be recovered as  $\CS = \Phi(\CP$).

b) In particular, 
every latitudinal gl star is topological.

c) Conversely, every symmetric gl pencil $\CP$ in $E$ yields a latitudinal gl star $\CS = \Phi(\CP$).
\ethm

\bpf
Assertions a) and c) are obvious. Note that the reflection of $E$ about $Z$ is induced by the element of order 2 in $\Phi$. 
Part b) follows from compactness of $\CP$ (Theorem \ref{pencil} (b)) and of $\Phi$, together with 
continuity of the surjective map $\Phi \times \CP \to \CS: (\phi,L) \to \phi(L)$.
\epf

The proof that latitudinal gl stars are topological given in \cite{latitud} uses four pages in print. This is because it 
consists of an explicit description of the map sending a point-line pair $(p,L)$ to the line parallel to $L$ and containing $p$. 
Verification of the continuity of several implicit maps is required in the course of this proof.\\


\section {The symmetric case}\label{sym case}

In order to specify a rotational gl star $\CS$ in standard position one has to enumerate the $\Phi$-orbits of lines in $\CS$.
The orbit consisting of $Z$ alone and the orbit consisting of the lines perpendicular to $Z$ and passing through the origin are 
always present. There may be other lines $L\in \CS$ meeting $Z$. In this case, the orbit of $L$ consists of all lines contained in some 
$\Phi$-invariant cone. If $L\in \CS$ misses $Z$, then its orbit is one of the two reguli contained in some 
hyperboloid invariant under $\Phi$. In order to distinguish between the two reguli, it is convenient to adopt this convention: 
We call $\Phi(L)$ a \it right regulus \rm if $L$ may be parametrized by a function $v: \BR \to L$ with $v(t) = (x(t),y(t), z(t))$ 
such that the lines $\BR (x(t),y(t))$ in the $(x,y)$-plane rotate counterclockwise and the function $z(t)$ is increasing.
If $z$ is decreasing, then we have a \it left regulus. \rm

Now $\CS$ may be described by exhibiting the points $\sigma(p_t)$ or the lines 
    $$L_t = p_t \vee \sigma(p_t) \quad {\rm for} \quad  p_t = (\sqrt{1-t^2},0,t), t \in [0,1].$$
A necessary condition is that $L_1 = Z$ and that $L_0$ passes through the origin and contains $p_0$; 
hence, $L_0$ is the $x$-axis $X$. Instead of exhibiting these lines, we may specify the cones and 
hyperboloids carrying their orbits, for example by writing down the hyperbolas or line pairs obtained by 
intersection with the plane $y=0$.  In this case, it is necessary to decide for each hyperpoloid whether the right or left 
regulus should be taken. In order to get a continuous involution $\sigma$ and a topological gl star, the choice of left or 
right regulus has to obey a rule: If  no lines $L_t$ meeting  $Z$ occur for $t_1 \le t \le t_2$, then 
the choice of left or right  
must not be changed within the interval  $[t_1,t_2]$.
In other words, a change from left to right or vice versa is only possible at parameters $t$ corresponding to a cone.  \\

A special case is that of \it symmetric rotational gl stars. \rm By symmetric, we mean symmetry 
of each of the cones and hyperboloids with respect to the $(x,y)$-plane. 
This means that the third coordinate of $\sigma(p_t)$ is always equal to $-t$, or that the line pairs or hyperbolas seen in the plane
$y = 0$ all have symmetry axes $Z$ and $X$. Note, however, that the gl star itself is usually \it not \rm symmetric with 
respect to the $(x,y)$-plane, because this symmetry map exchanges the right and left reguli in each hyperboloid. Instead, 
we have symmetry of the gl star about the $y$-axis (or in fact about any line passing through 0 and orthogonal to $Z$).

We remark here that the full automorphism group of the parallelism from a 
symmetric and rotational gl star is a degree four extension of $\BT^2$ and is isomorphic to the corresponding group in the case of a
latitudinal gl star. However, the actions are different. In the symmetric case, we have the reflection about the $y$-axis, 
while in the latitudinal case there is  the reflection about any plane containing $Z$, see \ref{char ax}. 
More details on full group actions are given in \cite{Stevin}. Note also that only 
Clifford paralllelism is both latitudinal and symmetric.\\

In the symmetric case, it is possible to describe all examples in an easy manner. \\

Let $a: [0,1[\ \to \ [0,\infty[$ be an increasing  bijection satisfying the conditions
\begin{equation} \label{eins}   
   t^2 \le \frac {a(t)^2}{1+ a(t)^2} \quad {\rm for \quad all} \quad t \in [0,1[, \quad {\rm and} 
\end{equation}
\begin{equation}\label{zwei}
   \lim_{0< t \to 0} \frac {t^2(1 + a(t)^2)}{a(t)^2} = 1 . 
\end{equation}   
Then 
    $$c(t)^2 := a(t)^2 - t^2(1+a(t)^2)$$ 
is nonnegative, and defines $c(t) \ge 0$. For $t > 0$, let $H_t$ be the hyperbola or line pair defined by the equation 
   $$a(t)^2x^2 - z^2 = c(t)^2.$$   
This yields a family of $\Phi$-invariant cones and hyperboloids.   

\bthm \label{sym}
\cite{Stevin} For any admissible choice of reguli, the cones and hyperboloids defined above yield a topological, symmetric and 
rotational gl star, and all such gl stars in standard position are obtained in this way. 
\ethm

\bpf
The hyperbola or line pair $H_t$ passes through the points $p_t^\pm = (\sqrt{1-t^2},\pm t)$ on the right half of the circle 
$\BS_1$ and intersects the positive $x$-axis at $v_t = (\frac {c(t)} {a(t)}, 0)$. The asymptotic slopes are $\pm a(t)$, and by
condition (\ref{zwei}), the points $v_t$ converge to the origin for $t \to 0$. 
This shows that the lines
$L_t$ contained in $\Phi (H_t)$ and passing through  $p_t^+$ converge to $X$ or $Z$ for $t \to 0$ or $t \to 1$, respectively. 
If any two lines in the resulting set $\CS$ meet outside the unit sphere, then there is a pair of hyperbolas $H_{t_1}$, $H_{t_2}$ 
with $t_1 < t_2$ whose upper right branches meet outside the circle. Now the asymptotes of these branches have slopes 
$a(t_1) < a(t_2)$, and the upper right intersection point of $H_{t_1}$ with the circle is below the corresponding point 
of the other hyperbola. This shows that the right upper branches of the two hyperbolas meet twice outside the circle, 
where we count a point of touching as two points of intersection. 
Now the hyperbolas are symmetric about both $X$ and $Z$, hence they have 8 points in common and coincide, a contradiction.
Now Theorem \ref{crit} implies that we have obtained a symmetric rotational gl star.

The necessity of the conditions imposed on the function $a$ are easily checked, and this proves the last statement of the theorem.
\epf


\section{The general case}\label{general}

The special gl stars of the previous two sections correspond to some set of functions between real intervals. 
For rotational gl stars in general, one expects that pairs of such functions are needed, because the 
points $\sigma(p_t)$ from the introduction to Section \ref {sym case} have to be specified.
In fact, Betten and Riesinger claim a result concerning a family of 
rotational gl-stars corresponding to pairs of such functions (\cite{Stevin}, Theorem 42), inclucing examples 
that are neither symmetric nor axial.
However, their proof is invalid because it uses the 
argument about intersecting hyperbolas from the  proof of Theorem \ref{sym} in a situation where the two hyperbolas 
are \it not \rm known to share two axes of symmetry. Here we give a proof which has virtually no overlap with  their attempt. 
Also we have corrected a rather obvious omission (surjectivity of the function $f$), and we admit both left and right 
reguli by introducing a sign function $\epsilon$. \\

Let $f : [0,1] \to [0,1]$ be an increasing (continuous) bijection and let $g: [0,1] \to [-1,0]$ 
be a continuous non-decreasing function 
satisfying  $g(0) = -1$, $g(1) = 0$ and the inequalities 
   $$-\sqrt{1-f(t)^2} \le g(t) \le 0  $$for all $t$. Moreover let $\epsilon: [0,1] \to \{1,-1\}$ be any function that is constant on every interval 
$I \subseteq [0,1]$ such that $-\sqrt{1- f(t)^2} < g(t)$ holds on $I$.
For the points  $p_t = (\sqrt{1-t^2}, 0, t)$, $t\in [0,1]$,    define
    $$\sigma_{f,g}(p_t) = (g(t), \epsilon(t)\sqrt{1- f(t)^2-g(t)^2}, -f(t)),$$
and complete $\sigma$ to a $\Phi$-invariant involution of $\BS_2$ in the unique way. 
Let $\CS_{f,g}$ be the resulting set of lines,
   $$\CS_{f,g} = \{u \vee \sigma_{f,g}(u) \mid u \in \BS_2\}.$$
   
\bthm \label{42}
The set $\CS_{f,g}$ defined above is a topological, rotational gl star.
\ethm
   
\bpf
1. Note that $\sigma$ is surjective and in fact an involution of $\BS_2$ because we have assumed that $f$ is surjective. 
This condition is missing in \cite {Stevin}. Note also that $\CS  = \CS_{f,g}$ contains the lines $L_t = p_t \vee \sigma(p_t)$ and,
in particular, contains the $z$-axis $Z = L_1$, the 
$x$-axis $L_0$, and the horizontal lines
$\phi(L_0)$, $\phi \in \Phi$. In particular, on the equator $\BS_2 \cap Z^\perp$ the involution 
$\sigma$ induces the antipodal map $p \to -p$.

2. We want to apply Proposition \ref{crit}. First we need to check that $\sigma$ is continuous. The function $\epsilon$ is 
constant except at zeros of $\sqrt{1- f(z)^2-g(t)^2}$, hence $\sigma$ is continuous on the 
half meridian $M = \{p_t \mid t \in [0,1]\}$ and on $\sigma(M)$, and
hence everywhere. 

3. By definition, no line in $\CS$ meets either the Z-axis or one of the horizontal lines $\phi(L_0)$ outside the sphere. 
So we only have to deal with exterior points of intersection between two of the remaining lines. By rotation invariance, 
we may always assume that one of the two lines is among the lines $L_t$. 

4. The sign factor $\epsilon$ serves to choose between left and right reguli. Since swapping left and right reguli 
does not affect the properties of a gl star, we may assume that $\epsilon$ is constant, say $\epsilon = -1$.\

5. If $f = id$, then $\CS$ is a symmetric rotational gl star by Theorem \ref{sym}, and if $g(t) = -\sqrt{1-f(t)^2}$ 
holds for all $t$, 
then  $\CS$ is an axial gl star by Theorem \ref{latitud}. Now suppose that the lines $L_t$ and $\phi(L_s)$ intersect outside $\BS_2$.
Then we can choose a number $r > \max\{s,t\}$ and change the functions $f$ and $g$ on the interval $[r,1]$ such that the condition 
for symmetric or axial gl stars is  satisfied on some subinterval $[u,1]$ with $u>r$. The new line set $\CS'$ still contains the same 
pair of intersecting lines. In order to prove that such a pair cannot occur, we may 
therefore assume that the local intersection condition of Theorem \ref{crit} is satisfied near the $z$-axis $Z$.

6. We proceed to show that the local intersection condition is satisfied near every line $L_t$, $t>0$. By rotational symmetry, 
it suffices to show that $L_t$ does not meet any line $\phi L_s$ outside the sphere, where $\phi \in \Phi$ rotates 
through an angle less than $\pi /2$ in positive or negative direction. Again by symmetry, we may also assume that $t < s$. 
Suppose that these lines meet in a point $v$ outside the sphere. We shall discuss the case that $v$ lies in the upper half space
defined by $z \ge 0$. The case $z \le 0$ is not identical, but sufficiently similar to 
allow for an analogous argument.

\begin{figure}[ht]
\begin{center}
\includegraphics[height=8cm]{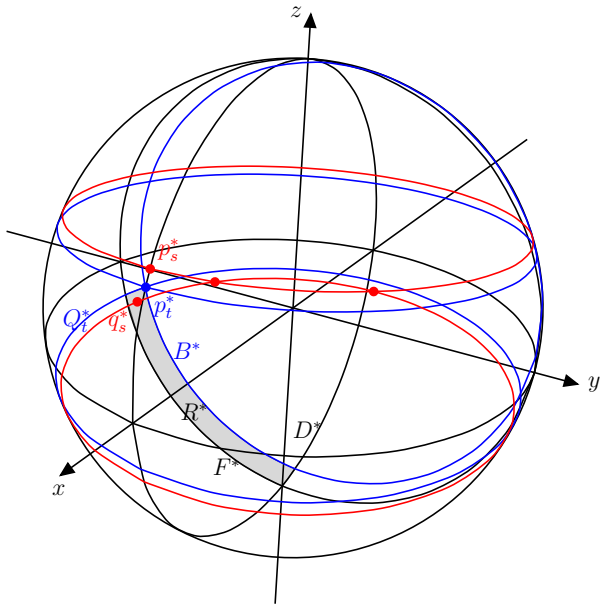}
\end{center}
\caption{Central projection from $v$ onto the plane $v^\perp$}\label{centralproj}
\end{figure}

7. We now look at a central projection of the situation, using $v$ as projection center and a suitable image plane $E$,
see Figure \ref{centralproj}.
For a point $p$, the projection is given by  $p^* = (p \vee v) \wedge E$. The fact that $v \in L_t$ means that 
$p_t^* = q_t^*$. Similarly, $v \in \phi L_s$ means that 
\begin{equation}\label {*}
  (\phi p_s)^* = (\phi q_s)^*. 
\end{equation}  
Now by the conditions on $f$ and $g$, 
we know that $q_s$ lies in a region $R$ on the sphere bounded above by the orbit circle $Q_t = \Phi q_t$, bounded backwards by the 
circle $B: x = g(t)$ and in front by the circle $F: x = 0$, and bounded below by the circle $D: y = 0$. Thus $q_s^*$ lies 
in the region $R^*$ bounded by the conics $Q_t^*, B^*, F^*$ and $D^*$. In particular, $q_s$ lies to the left of $B^*$. 
From equation (\ref{*}) we see that the orbits $P_s = \Phi p_s$ and $Q_s$ have intersecting images $P_s^*$ and $Q_s^*$. 
Then the orbits
$P_t$ and $Q_t$, which lie between the former two orbits, have  images intersecting in two points. These are $p_t^*$ 
and one further point,  situated to the right of $p_t^*$ on the lower arc of $P_t^*$. Since $\phi$ rotates through an 
angle less than $\pi/2$,
we see that $\phi$ is a clockwise rotation, so that $q_s^*$ is moved to the right. Then $p_s^*$ is moved to the left, and cannot
hit one of the intersection points of $P_s^*$ and $Q_s^*$, which lie between those of $P_t^*$ and $Q_t^*$. 
This shows that an intersection of $L_t$ and $\phi L_s$ in the upper half space is not possible. 

8. Neihgborhoods of $L_0$ are special. Such a neighborhood contains lines $\phi L_s$ with $s \ge 0$, where $\phi$ either rotates
through an angle near 0 or through an angle near $\pi$. Lines of the first kind obviously do not meet lines of the second kind 
outside the sphere. Lines of the same kind do not meet outside the sphere by what we have shown in the preceding steps.  
This completes the proof.
\epf

\bf Remark. \rm  The class of gl stars obtained in Theorem \ref{42} is substantial but does not comprise all rotational gl stars.
The restrictions arise from the monotonicity condition imposed on $g$, which is needed to make this proof work. Together with the 
necessary condition $g(1) = 0$, it enforces that $g$ is always negative and thus excludes many reguli from being part 
of the gl stars constructed here. However, we shall show later that every regulus consisting of 2-secants of the sphere can be part 
of a rotational gl star, see Theorem \ref{parabolas} below. \\

We return to the method of describing a topological rotational gl star presented at the beginning of Section \ref{sym case}. 
One starts from a set $\CH$ of cones and hyperboloids, each of them invariant under the action of the rotation group $\Phi$, 
and the gl star $\CS$ 
then consists of the axis $Z$, all horizontal lines containing the origin, and 
all lines contained in one of the cones, plus one of the two reguli contained in every one of the hyperboloids. 

\bthm \label{eqn-star}
In the $(x,z)$-plane $\BR \times 0 \times \BR$, consider a family of subsets $H_a$ defined by equations
\begin{equation} \label{hyp}
a^2x^2-(z-b)^2=c^2,
\end{equation}
for all $a >0$, where $b = b(a)$ and $0 \le c = c(a)$ are continuous functions of $a$. Assume that 
\begin{itemize}
   \item[(1)] 
    $b(a)^2 + c(a)^2 < a^2$ for all $a$.
   \item[(2)]
   $\lim_{a \to 0}b(a) = 0 = \lim _{a\to 0}\frac {c(a)} a $.   
   \item[(3)] 
    Every point $p_t = (\sqrt{1-t^2},0,t)$, $0\ne t \in \ ]-1,1[$, belongs to exactly one set $H_a$. 
   \item[(4)] 
   Every point $(x,0,z)$ with $x > 0$, $z \ne 0$ and $x^2 + z^2 \ge 1$ belongs to at most one set $H_a$.
\end{itemize}
Then we obtain a topological, rotational gl star $\CS$, consisting of the $z$-axis $Z$, the orbit $\Phi X$ of the $x$-axis $X$, and 
all lines contained in cones $\Phi H_a$, $c(a)=0$, plus the right reguli contained in the hyperboloids $\Phi H_a$ for $c(a)\ne 0$. 
On intervals bounded by consecutive zeros of the function $c$, right reguli may be replaced by left reguli. 

Conversely, every topological, rotational gl star in standard position is obtained in this way. 
\ethm

\bf Remark. \rm The axial and symmetric cases are given by $c \equiv 0$ and by $b \equiv 0$, respectively. 
It is not hard to recover the earlier
Theorems \ref{latitud} and \ref {sym} from the present result in those special cases.

\bpf
1. Condition (1) says that every set $H_a$ intersects $X$ inside the unit circle. Consequently, $H_a$ 
contains a point $p_{t(a)}$ with $t(a) > 0$.
The function $t: \ ]0,\infty[ \ \to \ ]0,1[$ is continuous and, by (3), bijective. Condition (4) 
implies that $t$ is not decreasing, so it is an 
increasing function of $a$. Every point $p_t$, $t \in  \ ]0,1[$, satisfies a unique equation (\ref{hyp}), 
whose coefficients are continuous functions of $t$.
There is a unique line $L_t \in \CS$ which satisfies the same equation and belongs to the right 
regulus of $\Phi H_a$ if $c \ne 0$. Therefore, the map $L: t \to L_t$ is continuous. 
By (1), $L_t$ meets the interior of the sphere $\BS_2$, hence $L_t$ intersects the sphere in a second point $q_t = :\sigma_1(p_t)$. 

In step 2, we shall show that the map $L :\ ]0,1[ \to \CS$ is continuously extended to $[0, 1]$ by setting $L_0 = X$ and $L_1 = Z$. 
Then it will follow that $\sigma_1$ extends uniquely to a $\Phi$-invariant map $\sigma_2: \BS_2^+ \to \BS_2^-$ 
from the (closed) upper hemisphere to the lower hemisphere. By (4) and by $\Phi$-invariance, $\sigma_2$ is injective. 
On the equator $\BS_2^+ \cap \BS _2^-$, $\sigma_2$ 
induces the antipodal map 
$-id$. It follows that $\sigma_2$ is surjective, because the antipodal map is not null homotopic in any proper 
subset of the lower hemisphere, or alternatively, because the connected image of $\sigma_2$ is $\Phi$-invariant and contains 
both the equator and $\sigma_2(p_1) = -p_1$. Thus we obtain a 
fixed point free involution $\sigma: \BS_2 \to \BS_2$ extending $\sigma_2$, such that  $\CS$ consists of the 
lines $p \vee \sigma (p)$, $p\in \BS_2$. 

In order to complete the proof by applying
Proposition \ref{crit}, we need to know that two lines in $\CS$  never meet outside the sphere. This follows from condition (4)
if the lines are contained in two distinct sets $\Phi H_a$. It follows from (1) if one of the lines belongs to $\Phi X$. Lines of the 
same regulus are always disjoint. The last remaining case arises for cones $H_a$. Two lines of the cone intersect in the 
vertex, which  moreover belongs to $Z \in \CS$. So we have to show that all vertices of cones are inside the sphere. This follows from
(4) together with our claim that $L_t \to Z$ if $t \to 1$. 

2. If $t \to 1$ then $p_t \to p_1 \in Z$. The line $L_t$ contains $p_t$. Since $a(t) \to \infty$, the point at infinity 
of $L_t$ belongs to a $\Phi$-orbit that converges to the point at infinity of $Z$. This proves that $L_t \to Z$. 
For $t \to 0$, we have that $p_t \to p_0 \in X$. Moreover, if $c = 0$, then $L_t$ contains the vertex of the cone containing 
$p_t$, and the vertex converges to the orgin since $b \to 0$. For those $t$ with $c \ne 0$, the line 
$L_t$ contains a point from the $\Phi$-orbit of the vertex $v= (b,\frac c a )$ of the hyperbola $H_a$. Hence,
the distance of $L_t$ from the origin is at most $\Vert v \Vert$, and $L_t \to X$ follows from (2).

3. Conversely, suppose that we ave a topological rotational gl star $\CS$. Then we obtain a system of hyperboloids and cones and a family of equations (\ref{hyp}).
The necessity of conditions (3) and (4) is obvious, as they are part of the definition of a gl star. 
Condition (1) expresses the fact that $H_a$ 
intersects $X$ inside the unit circle, which follows from $X \in \CS$. Condition (2) is obtained by reversing the 
argument in step 2. This proves the last assertion.
\epf

We wish to reshape the preceding result, describing the sets $H_a$ in terms of their intersection points with the unit circle rather
than by the coefficients of their equation. This  might help to guess examples. 

Consider the situation of Theorem \ref {eqn-star} and retain the notation used above. The set $H_a$ intersects the unit circle
in 4  points. Those with a positive $x$-coordinate are $p_{t(a)}$ and $p_{-s(a)}$, where 
\begin{align} 
 t(a) &=\frac{b+\sqrt{b^2+(a^2+1)(a^2-b^2-c^2)}}{a^2+1}, \\
 -s(a) &=\frac{b-\sqrt{b^2+(a^2+1)(a^2-b^2-c^2)}}{a^2+1}.
 \end{align}
As before, the functions $t$ and $s$ are homeomorphisms $]0,\infty[ \ \to \ ]0,1[$. The coefficients $b$ and $c$ are 
expressed in terms of $t$ and $s$ by the equations
\begin{align}
b(a) &=(a^2+1)\left(\frac{t-s}{2}\right),\label{eq:b}\\
c(a)^2&=a^2-(a^2+1)\left(\left(\frac{t+s}{2}\right)^2+a^2\left(\frac{t-s}{2}\right)^2\right).\label{eq:c}
\end{align}
This makes sense if the expression for $c^2$ is nonnegative. Thus, 
the functions $t$ and $s$ are restricted by the condition
\begin{align*}
&a^2-(a^2+1)\left(\left(\frac{t+s}{2}\right)^2+a^2\left(\frac{t-s}{2}\right)^2\right)\\ 
&=-\frac14\left((t-s)^2a^4-2(2-t^2-s^2)a^2+(t+s)^2\right)\\
&\ge 0
\end{align*}
for all $a>0$. This inequality can be rewritten as 
\begin{equation}
\frac{a^2}{a^2+1}-ts\ge (a^2+1)\left(\frac{t-s}{2}\right)^2.\label{c2s}
\end{equation}
With the formulas (\ref{eq:b}) and (\ref{eq:c}) for $b$ and $c$ the equation for $H_a$ becomes
\begin{equation}
h_{x,z}(a):=(a^2+1)(x^2-1+(t-s)z+ts)+1-x^2-z^2=0.\label{has}
\end{equation}
Using these results, we may rewrite Theorem \ref{eqn-star} as follows.

\bcor\label{param-star}
Let $t,s:[0,\infty[\ \to [0,1[$ be homeomorphisms such that 
\begin{itemize}
	\item[(1)]
	$\lim_{0 \ne a\to 0} \frac{t(a) + s(a)}{2a} = 1$.
	\item[(2)]
	$\frac{a^2}{a^2+1}-t(a)s(a)\ge (a^2+1)\left(\frac{t(a)-s(a)}{2}\right)^2
	$ for all $a>0$.
	\item[(3)]
	The map $h_{x,z}$ given by
	\begin{align*}
a\mapsto &(a^2+1)(x^2-1+(t(a)-s(a))z+t(a)s(a))+1-x^2-z^2\\
&=a^2(x^2+z^2-1)+(a^2+1)(t(a)-z)(s(a)+z)
\end{align*} 
has at most one positive root for each $(x,z)$ such that $x>0$, $z\ne 0$ and $x^2+z^2 \ge 1$. 
\end{itemize}
Then we obtain a topological, rotational gl star by applying Theorem \ref{eqn-star} to the values  $b(a)$ and $c(a)$ 
defined by equations (\ref {eq:b}) and (\ref{eq:c}) above.

Conversely, every topological, rotational gl-star in standard position can be obtained in this way for 
suitable functions $t$ and $s$ as above.
\ecor

\bf Remark. \rm Here, the symmetric case is given by $s = t$ and  the axial case by $c \equiv 0$. 

\bpf
Note first that we have required continuity of $t(a)$ and $s(a)$ at $a=0$, whereas continuity of $b(a)$ and $c(a)$ in Theorem 
\ref{eqn-star} is required only for $a > 0$. Therefore, $\lim_{a \to 0} b = 0$ holds. Moreover, 
$\lim_{0\ne a\to 0} \frac c a = 0$ follows easily from  the present condition (1). From equations 
(\ref {eq:b}) and (\ref{eq:c}) we get $0 \le b^2+c^2=a^2-(a^2+1)ts \le a^2$, which gives condition (1) of Theorem \ref{eqn-star}.
Conditions (4) of the theorem is  condition (3) of the corollary, and condition (3) of the theorem is guaranteed because
the surjective functions $t$ and $s$ describe, by their definition, the $z$-coordinates of the intersections of $H_a$ 
with the unit circle.
\epf

\bf Example. \rm Checking the properties of the functions $t$ and $s$ can be very hard. We succeded in one concrete example.
For each $r>0$ let $\varphi_r:[0,\infty[ \ \to [0,1[$ be the map defined by
$$\varphi_r(a)=\frac{a(a+r)}{a^2+ra+r}=1-\frac{r}{a^2+ra+r}.$$
If we set 
$$t=\varphi_{\frac 3 2}\quad \text{and}\quad s=\varphi_2,$$
then after some lengthy computation we obtain
$$h_{x,z}(a)=\frac{l(a)}{(2 a^2+3 a+3) (a^2+2 a+2)},$$
where 
\begin{align*}
l(a)&=2x^2a^6 +7x^2a^5 +(13 x^2-2 z^2
+ z-5)a^4\\
&\qquad +(12 x^2-7 z^2-5)a^3 +(6 x^2-13 z^2+z)a^2 
-12 a z^2-6 z^2,
\end{align*}
Using Descartes' rule of signs, it can be shown that this function has exactly one positive root for each 
admissible pair $(x,z)$. The other conditions of the Corollary can  be checked without difficulty and one 
gets a rotational gl star. 

By continuity, the choice $t = \phi_r$ and $s = \phi_{r'}$, where $(r,r')$ is sufficiently close to $(\frac 3 2 ,2)$ also 
yields a rotational gl-star.\\

\bf Remark. \rm
Here is a different perspective on condition (3) of the corollary. We focus attention on a horizontal line $z = z_0 \ne 0$. 
If we let
$$a_z=\begin{cases} \infty, &\text{if }|z|\ge 1,\\
t^{-1}(z), &\text{if }0<z<1,\\
s^{-1}(-z), &\text{if }-1<z<0,
\end{cases}$$
and
$$p_z(a)=\frac{a^2+1}{a^2}(z-t(a))(z+s(a))$$
where $z\ne 0$, then condition (3) from Corollary \ref{param-star} is equivalent to $p_z$ being injective (in fact, 
strictly decreasing) on the interval $]0,a_z[$ for all $z\ne 0$. 

Indeed, $h_{x,z}(a)=0$ yields $x^2+z^2-1=\frac{a^2+1}{a^2}(z-t(a))(z+s(a))$, and every $(x,z)$ with $x>0$ and $x^2>1-z^2$ 
is covered  by the sets $H_a$ at most once. Only those $a$ for which $(z-t(a))(z+s(a))>0$ 
are relevant, which leads to the definition of $a_z$. \\

We conclude with a fairly versatile construction based on Theorem \ref {eqn-star}. We shall obtain a large set of examples which 
shows, in particular, that the family of examples provided by Theorem \ref{42} is far from exhausting all possibilities.

We shall use a technique of linear interpolation which is awkward to handle for hyperbolas. Therefore we use a transformation which 
translates families of hyperbolas into families of parabolas. As always, we allow degenerate hyperbolas (i.e., line pairs). 
So let $\omega : \BR \times 0 \times \BR \to \BR^2$ be given by
   $$ (u,v) = \omega(x,0,z) = (z, x^2).$$
This maps the  hyperbola  (or line pair) $H$ with equation
       $$ H: \quad a^2x^2 - (z-b)^2 = c^2$$
(where $0\le c$ and $0 <a$) to the parabola $P = \omega H$ with equation
      $$P: \quad v = \alpha (u - \beta)^2 + \gamma,$$
where
       $$\alpha = \frac 1 {a^2} \quad \beta = b, \quad \gamma = \frac {c^2}{a^2}.$$   
Moreover, $\omega$ sends the vertices $(b,\pm \frac c a)$ of $H$ to the vertex $(\beta,\gamma)$  of $\omega H$.     
The $z$-axis $Z$ is mapped to the $u$-axis $U$, and the right half ($x \ge 0$) of the unit circle $C$ goes to the  parabolic arc 
$D=\{(u,1-u^2)\mid -1\le u\le 1\}$. Let $R$ be the bounded region bounded by  $U \cup D$. 

\begin{figure}[ht]
\begin{center}
\includegraphics[height=4cm]{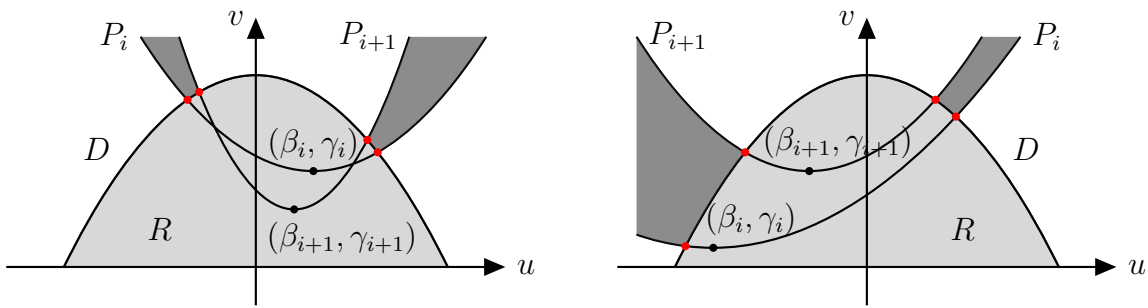}
\end{center}
\caption{The parabolas $P_i$ and $P_{i+1}$}\label{fig:parabolas}
\end{figure}

We choose a two-way infinite sequence of parabolas $P_i$ indexed by integers $i$, defined by  equations 
$v = p_i(u) = \alpha_i (u - \beta_i)^2 + \gamma_i$ with $\alpha_i>0$ and $\gamma_i\ge 0$, 
such that the following conditions hold, compare Figure \ref{fig:parabolas}:
\begin {itemize}
\item[(1)]
The sequence $\{a_i\}_i$ is strictly increasing with $\lim_{i \to -\infty}a_i = 0$ and $\lim_{i\to\infty}a_i = \infty$.
\item[(2)]
$\lim_{i \to -\infty} \gamma_i = 0.$
\item[(3)]
$P_i \cap D$ consists of two points separated by the $v$-axis $V$, and each point of $P_{i+1} \cap D$ separates the 
points in $P_i \cap D$ on the arc $D$.
\item[(4)]
$P_i$ and $P_{i+1}$ do not intersect outside $R$.
\item[(5)]
The vertices $(\beta_i, \gamma_i)$ of $P_i$ converge to the origin as $i \to \infty$. 
\end{itemize}

Now we interpolate $i \to p_i$ in order to obtain a one-parameter family $P_r, r \in \BR$, of parabolas. 
That is, for $t \in [0,1]$ we define $P_{i+t}$ by
    $$p_{i+t} : = (1-t) p_i +  tp_{i+1}.$$
These parabolas simply cover the two regions of $\BR^2 \setminus \overline R$ bounded by $P_i$ and $P_{i+1}$, and the 
limit conditions (1) and (3) now hold for $r \to \pm \infty$.  This is obvious for (1), and for (3) we use
the following  observation.  For $u < \min\{\beta_i,\beta_{i+1}\}$, the functions $p_i$ and $p_{i+1}$ are both decreasing, and the 
same holds for $p_{i+t}$, $t \in  [0,1]$. Thus the vertex $(u_t,v_t)$ of $p_{i+t}$ satisfies $u_t \ge \min \{\beta_i,\beta_{i+1}\}$,
and likewise, $u_t \le \max \{\beta_i,\beta_{i+1}\}$. Moreover, the function $p_{i+t}$ lies between the functions $p_i$ and $p_{i+1}$,
hence its minimum value is bounded above by the larger one of their minimum values.
This proves our claim, and we have the following.

\bthm \label{parabolas}
Given a sequence  $P_i$ of parabolas satisfying conditions (1), (2) and (3) as 
above, taking $\omega$-inverse images of the interpolated 
family $P_r$, $r \in \BR$, we obtain a family of hyperbolas or line pairs as in Theorem \ref{eqn-star}, and thus we obtain a 
topological, rotational gl star.
\ethm

\bpf
Remember that $\alpha = \frac 1 {a^2}$ is a strictly increasing function of $r \in \BR$, hence the family of 
hyperbolas may be parametrized by $a > 0$.
Only condition (1) of \ref{eqn-star} then needs an explanation: this condition expresses that the hyperbolas intersect the axis $X$
inside the unit circle. Under $\omega$, this corresponds to the property that the parabolas intersect the $v$-axis $V$ inside the 
region $R$ bounded by the parabolic curve $D$ and the axis $U$. This is guaranteed by our condition (2) above. 
\epf

\bf Remarks. \rm 1. As announced earlier, this result shows that Theorem \ref{42} is far from giving a complete view 
of all rotational gl stars in standard position. We see here that every hyperboloid meeting the $(x,y)$-plane inside the unit sphere 
carries two reguli that may be embedded in rotational gl stars. By contrast, in Theorem \ref{42} the hyperboloids 
are restricted by the incisive condition $g(t) \le 0$. 

Moreover, it is easy to choose the parabolas such that the vertices 
of the corresponding hyperbolas lie outside the unit circle for large values of the slope $a$. They may even tend to 
$\infty$ as $a \to \infty$.

2. It is easy to modify the construction so that one may start with a finite sequence or a one-way infinite sequence of parabolas. 
One just completes the family by adding $tp_i$, $t \in [0,1]$, if the sequence starts with index $i$, and by 
adding $tp_j$, $t \ge 1$, if the sequence ends with index $j$ and the vertex of $P_j$ is the origin.


\bibliographystyle{plain}

\bigskip
\bigskip
\noindent{Rainer L\"owen\\ 
Institut f\"ur Analysis und Algebra\\
Technische Universit\"at Braunschweig\\
Pockelsstr. 14\\
38106 Braunschweig\\
Germany\\ \\
G\"unter F. Steinke\\
           School of Mathematics and Statistics\\
           University of Canterbury\\
	   Private Bag 4800\\
           Christchurch 8140, New Zealand}

\end{document}